\documentclass[leqno]{article}
\usepackage{a4wide,amsmath,amssymb,amsfonts,exscale}
\long\def\beginFORGET#1\endFORGET{}

\def\Ksep{{K^{\rm sep}}}

\def\an{{\rm an}}

 \def\Gal{\mathop{\rm Gal}\nolimits}
 \def\End{\mathop{\rm End}\nolimits}

 \def\Spec{\mathop{\rm Spec}\nolimits}

 \def\mod{\mathop{\rm mod}\nolimits}

\def\Stab{\mathop{\rm Stab}\nolimits}

\def\Ad{\mathop{\rm Ad}\nolimits}

\def\trace{\mathop{\rm tr}\nolimits}

\def\GL{\mathop{\rm GL}\nolimits}
\def\SL{\mathop{\rm SL}\nolimits}

\def\PGL{\mathop{\rm PGL}\nolimits}

\def\geom{{\rm geom}}

\def\sep{{\rm sep}}

\def\id{{\rm id}}

\let\phi\varphi
\let\epsilon\varepsilon

\newtheorem{Thm}{Theorem}[section]
\newtheorem{Prop}[Thm]{Proposition}
\newtheorem{Lem}[Thm]{Lemma}

\newtheorem{Def}[Thm]{Definition}

\newtheorem{Ques}[Thm]{Question}

\def\qed{{\hskip0pt\unskip\unskip\nobreak\hfil\penalty50
          \hskip1em\hbox{}\nobreak\hfil
           {$\square$}
          \parfillskip=0pt\finalhyphendemerits=0
          \par}\medskip}

\newenvironment{Proof}
               {\noindent{\bf Proof.}\ }
               {\qed}

\newenvironment{Proofof}[1]
               {\noindent{\bf Proof of #1.}\ }
               {\qed}

\newcommand{\BA}{{\mathbb{A}}}

\newcommand{\BC}{{\mathbb{C}}}

\newcommand{\BF}{{\mathbb{F}}}

\newcommand{\BP}{{\mathbb{P}}}

\newcommand{\Fa}{{\mathfrak{a}}}

\newcommand{\Fp}{{\mathfrak{p}}}

\newbox\mybox
\def\arrover#1{\mathrel{
       \setbox\mybox=\hbox spread 1.4em
              {\hfil$\scriptstyle#1$\hfil}
       \vbox{\offinterlineskip\copy\mybox
             \hbox to\wd\mybox{\rightarrowfill}}}}

\def\larrover#1{\mathrel{
       \setbox\mybox=\hbox spread 1.4em
              {\hfil$\scriptstyle#1\vphantom{g}$\hfil}
       \vbox{\offinterlineskip\copy\mybox
             \hbox to\wd\mybox{\leftarrowfill}}}}

\def\ontoover#1{\mathrel{
       \setbox\mybox=\hbox spread 1.4em
              {\hfil$\scriptstyle#1\vphantom{g}$\hfil}
       \vbox{\offinterlineskip\copy\mybox
             \hbox to\wd\mybox{\rightarrowfill\hskip-2.8mm
                               $\rightarrow$}}}}
\def\leftontoover#1{\mathrel{
       \setbox\mybox=\hbox spread 1.4em
              {\hfil$\scriptstyle#1\vphantom{g}$\hfil}
       \vbox{\offinterlineskip\copy\mybox
             \hbox to\wd\mybox{$\leftarrow$\hskip-2.8mm
                               \leftarrowfill}}}}
\let\longto\longrightarrow
\let\into\hookrightarrow
\let\onto\twoheadrightarrow

\def\isoto{\arrover{\sim}}

\begin{document}


\title{Monodromy Groups associated to Non-Isotrivial \\
Drinfeld Modules in Generic Characteristic}

\author{Florian Breuer\footnotemark[1] \and Richard Pink\footnotemark[7]}

\date{\today}

\maketitle

\renewcommand{\thefootnote}{\fnsymbol{footnote}}
\footnotetext[1]{Dept. of Mathematics, University of Stellenbosch, 
Stellenbosch 7600, South Africa, flo@math.jussieu.fr}
\footnotetext[7]{Dept. of Mathematics, ETH Zentrum, 8092 Z\"urich, 
Switzerland, pink@math.ethz.ch}

\begin{abstract}
Let $\phi$ be a non-isotrivial family of Drinfeld $A$-modules of rank $r$ 
in generic characteristic
with a suitable level structure over a connected smooth algebraic variety $X$. 
Suppose that the endomorphism ring of $\phi$ is equal to~$A$. 
Then we show that the closure of the analytic fundamental group
of $X$ in $\SL_r(\BA_F^f)$ is open, where $\BA_F^f$ denotes the ring of 
finite ad\`eles of the quotient field $F$ of~$A$.

From this we deduce two further results: 
(1) If $X$ is defined over a finitely generated field extension of~$F$, 
the image of the arithmetic \'etale fundamental group of $X$ on the ad\`elic 
Tate module of $\phi$ is open in~$\GL_r(\BA_F^f)$.
(2) Let $\psi$ be a Drinfeld $A$-module of rank $r$ defined over a finitely 
generated field extension of~$F$, and suppose that $\psi$ cannot be defined 
over a finite extension of~$F$.
Suppose again that the endomorphism ring of $\psi$ is~$A$.
Then the image of the Galois representation on the ad\`elic Tate module of 
$\psi$ is open in $\GL_r(\BA_F^f)$. 

Finally, we extend the above results to the case
of arbitrary endomorphism rings.
\end{abstract}

\begin{flushleft}

\end{flushleft}

\parindent=0pt


\section{Analytic monodromy groups}
\label{1}


Let $\BF_p$ be the finite prime field with $p$ elements.  Let $F$ be a
finitely generated field of transcendence degree $1$ over~$\BF_p$. 
Let $A$ be the ring of elements of $F$ which are regular outside a
fixed place $\infty$ of~$F$.  
Let $M$ be the fine moduli space over $F$ of Drinfeld $A$-modules of rank $r$ with some sufficiently high level structure.
This is a smooth affine scheme of dimension $r-1$ over~$F$.

\medskip
Let $F_\infty$ denote the completion of $F$ at $\infty$, and $\BC$ the completion of an algebraic closure of~$F_\infty$.
Then the rigid analytic variety $M_\BC^\an$ is a finite disjoint union of
spaces of the form $\Delta\backslash\Omega$, where $\Omega \subset (\BP^{r-1}_\BC)^\an$ is
Drinfeld's upper half space and $\Delta$ is a congruence subgroup of $\SL_r(F)$
commensurable with $\SL_r(A)$.

\medskip
Let $X_\BC$ be a smooth irreducible locally closed algebraic subvariety of~$M_\BC$. 
Then $X_\BC^\an$ lies in one of the components $\Delta\backslash\Omega$ of~$M_\BC^\an$.
Fix an irreducible component  $\Xi\subset\Omega$ of the pre-image of~$X_\BC^\an$. Then $\Xi \to X_\BC^\an$ is an unramified Galois covering with Galois group $\Delta_\Xi := \Stab_\Delta(\Xi)$.

\medskip
Let $\phi$ denote the family of Drinfeld modules over~$X_\BC$ determined 
by the embedding $X_\BC\subset M_\BC$. 
We assume that $\dim X_\BC\ge1$. Since $M$ is a fine moduli space, this means that $\phi$ is non-isotrivial. 
It also implies that $r\ge2$.
Let $\eta_\BC$ be the generic point of $X_\BC$ and $\bar\eta_\BC$ a geometric point above it.
Let $\phi_{\bar\eta_\BC}$ denote the pullback of $\phi$ to~$\bar\eta_\BC$.
Let $\BA_F^f$ denote the ring of finite ad\`eles of~$F$.
The main result of this article is the following:

\begin{Thm}
\label{AnZD}
In the above situation,
if $\End_{\bar\eta_\BC}(\phi_{\bar\eta_\BC}) = A$, then the closure of $\Delta_\Xi$ in $\SL_r(\BA_F^f)$ is an open subgroup of $\SL_r(\BA_F^f)$.
\end{Thm}

The proof uses known results on the $\Fp$-adic Galois representations associated to Drinfeld modules \cite{PinkDMT} and on strong approximation \cite{PinkSA}.

\medskip

Theorem \ref{AnZD} leaves open the following natural question:

\begin{Ques}
\label{GammaArith}
If $\End_{\bar\eta_\BC}(\phi_{\bar\eta_\BC}) = A$, is $\Delta_\Xi$ an arithmetic subgroup of $\SL_r(F)$?
\end{Ques}




\medskip

Theorem \ref{AnZD} has applications to the analogue of 
the Andr\'e-Oort conjecture for Drinfeld moduli spaces: 
see~\cite{BreuerAO}. Consequences for \'etale monodromy groups and for Galois 
representations are explained in Sections \ref{2} and~\ref{3}.
The proof of Theorem~\ref{AnZD} will be given in Sections \ref{p} through~\ref{c}.
Finally, in Section \ref{End} we outline the case of arbitrary endomorphism rings.

\medskip
For any variety $Y$ over a field $k$ and any extension field $L$ of $k$ we will abbreviate $Y_L := Y\times_k L$.


\section{\'Etale monodromy groups}
\label{2}


We retain the notations from Section~\ref{1}. 
Let $k \subset \BC$ be a subfield that is finitely generated over~$F$, such that $X_\BC = X \times_{k} \BC$ for a subvariety $X \subset M_k$.
Let $K$ denote the function field of $X$ and $K^\sep$ a separable closure of~$K$.
Then $\eta:=\Spec K$ is the generic point of~$X$ and $\bar\eta := \Spec K^\sep$ a geometric point above~$\eta$.
Let $k^\sep$ be the separable closure of~$k$ in~$K^\sep$.
Then we have a short exact sequence of \'etale fundamental groups
$$1\longto 
\pi_1(X_{k^\sep},\bar\eta)
\longto \pi_1(X,\bar\eta) \longto 
\Gal(k^\sep/k) \to 1.$$
Let $\hat A \cong \prod_{\Fp\not=\infty} A_\Fp$ denote the profinite completion of~$A$. 
Recall that $\BA_F^f \cong F\otimes_A\hat A$ and contains $\hat A$ as an open subring.
Let $\phi_\eta$ denote the Drinfeld module over $K$ corresponding to~$\eta$.
Its ad\`elic Tate module $\hat T(\phi_\eta)$
is a free module of rank $r$ over~$\hat A$. Choose a basis and let
$$\rho:\, \pi_1(X,\bar\eta) \longto \GL_r(\hat A) \subset \GL_r(\BA_F^f)$$
denote the associated monodromy representation.
Let $\Gamma^\geom \subset \Gamma \subset \GL_r(\hat A)$ denote the images of
$\pi_1(X_{k^\sep},\bar\eta)
\subset \pi_1(X,\bar\eta)$ under~$\rho$.


\begin{Lem}
\label{claim0a}
$\Gamma^\geom$ is the closure of $g^{-1}\Delta_\Xi\, g$ in $\SL_r(\hat A)$
for some element $g \in \GL_r(\BA^f_F)$.
\end{Lem}

\begin{Proof}
Choose an embedding $K^\sep\into\BC$ and a point $\xi\in\Xi$ above~$\bar\eta$. 
Let $\Lambda\subset F^r$ be the lattice corresponding to the Drinfeld module at~$\xi$. 
This is a finitely generated projective $A$-module of rank~$r$. 
The choice of a basis of $\hat T(\phi_\eta)$ yields a composite embedding
$$\hat A^r \ \cong\ \hat T(\phi_\eta) \ \cong\ \Lambda\otimes_A\hat A
\ \into\ F^r\otimes_A\hat A \ \cong\ (\BA_F^f)^r,$$
which is given by left multiplication with some element  $g \in \GL_r(\BA^f_F)$.
Since the discrete group $\Delta \subset \SL_r(F)$ preserves~$\Lambda$, 
we have $g^{-1}\Delta g \subset \SL_r(\hat A)$. 

\medskip
For any non-zero ideal $\Fa\subset A$ let $M(\Fa)$ denote the moduli space obtained from $M$ by adjoining a full level $\Fa$ structure. 
Then $\pi_\Fa\!: M(\Fa) \onto M$  is an \'etale Galois covering with group 
contained in $\GL_r(A/\Fa)$, 
and one of the connected components of $M(\Fa)_\BC^\an$ above the connected component $\Delta\backslash\Omega$ of $M_\BC^\an$ has the form $\Delta(\Fa)\backslash\Omega$ for 
$$\Delta(\Fa) := \bigl\{ \delta\in\Delta 
               \bigm| g^{-1}\delta g \equiv \id\, \mod \Fa\hat A \bigr\}.$$
Let $X(\Fa)_{k^\sep}$ be any connected component of the inverse image $\pi_\Fa^{-1} (X_{k^\sep})\subset M(\Fa)_{k^\sep}$. Since $k^\sep$ is separably closed, the variety $X(\Fa)_\BC$ over $\BC$ obtained by base change
is again connected.
The associated rigid analytic variety $X(\Fa)_\BC^\an$ is then also connected 
(cf. \cite[Kor.$\,$3.5]{Luetkebohmert}) 
and therefore a connected component of $\pi_\Fa^{-1} (X_\BC^\an)$. But one of these connected components is $\bigl( \Delta_\Xi \cap \Delta(\Fa) \bigr) \backslash \Xi$, whose Galois group over $X_\BC^\an \cong \Delta_\Xi\backslash\Xi$ is $\Delta_\Xi / \bigl( \Delta_\Xi \cap \Delta(\Fa) \bigr)$. 
This implies that $g^{-1}\Delta_\Xi\, g$ and $\pi_1(X_{k^\sep},\bar\eta)$ have the same images in $\GL_r(A/\Fa) = \GL_r(\hat A/\Fa\hat A)$. By taking the inverse limit over the ideal $\Fa$ we deduce that the closure of $g^{-1}\Delta_\Xi\, g$ in $\SL_r(\hat A)$ is~$\Gamma^\geom$, as desired.
\end{Proof}


\begin{Lem}
\label{EE}
$\End_{K^\sep}(\phi_{\eta}) = \End_{\bar\eta_\BC}(\phi_{\bar\eta_\BC})$.
\end{Lem}

\begin{Proof}
By construction $\bar\eta_\BC$ is a geometric point above~$\eta$, and $\phi_{\bar\eta_\BC}$ is the pullback of~$\phi_\eta$.
Any embedding of $K^\sep$ into the residue field of $\bar\eta_\BC$ induces a morphism
$\bar\eta_\BC \to \bar\eta$.
Thus the assertion follows from the fact that for every Drinfeld module over a field, any endomorphism defined over any field extension is already defined over a finite separable extension.
\end{Proof}


\begin{Thm}
\label{Main1}
In the above situation, suppose that $\End_{K^\sep}(\phi_{\eta}) = A$. Then
\begin{itemize}
\item[(a)] 
$\Gamma^\geom$ is an open subgroup of $\SL_r(\BA_F^f)$, and
\item[(b)] $\Gamma$ is an open subgroup of $\GL_r(\BA_F^f)$.
\end{itemize}
\end{Thm}

\begin{Proof}
By Lemma \ref{EE} the assumption implies that $\End_{\bar\eta_\BC}(\phi_{\bar\eta_\BC}) =  A$.
Thus part (a) follows at once from Theorem~\ref{AnZD} and Lemma~\ref{claim0a}.
Part (b) follows from (a) and the fact that $\det(\Gamma)$ is open in 
$\GL_1(\BA_F^f)$.
This fact is a consequence of work of Drinfeld \cite[\S8 Thm.$\,$1]{Drinfeld}
and Hayes \cite[Thm.$\,$9.2]{HayesExplicit} on the abelian class field 
theory of~$F$,
and of Anderson \cite{Anderson} on the determinant Drinfeld module. 
Note that Anderson's paper only treats the case $A=\BF_q[T]$; the general
case has been worked out by van der Heiden \cite[Chap.$\,$4]{vdHeiden}.
Compare also \cite[Thm.$\,$1.8]{PinkDMT}.
\end{Proof}


\section{Galois groups}
\label{3}

Let $F$ and $A$ be as in Section~\ref{1}. 
Let $K$ be a finitely generated extension field of $F$ of arbitrary trans\-cen\-dence degree, 
and let $\psi:\, A\to K\{\tau\}$ be a Drinfeld $A$-module of rank $r$ over~$K$.
Let $K^\sep$ denote a separable closure of~$K$ and
$$\sigma:\, \Gal(\Ksep/K) \longto \GL_r(\BA_F^f)$$
the natural representation on the ad\`elic Tate module of~$\psi$.
Let $\Gamma \subset \GL_r(\BA_F^f)$ denote its image.

\begin{Thm}
\label{DMT2}
In the above situation, suppose that $\End_{K^\sep}(\psi) = A$
and that $\psi$ cannot be defined over a finite extension of $F$ inside $K^\sep$.
Then $\Gamma$ is an open subgroup of $\GL_r(\BA_F^f)$.
\end{Thm}

\begin{Proof}
The assertion is invariant under replacing $K$ by a finite extension. We may therefore assume that $\psi$ possesses a sufficiently high level structure over~$K$. Then $\psi$ corresponds to a $K$-valued point on the moduli space $M$ from Section~\ref{1}. Let $\eta$ denote the underlying point on the scheme~$M$, and let $L\subset K$ be its residue field. Then $\psi$ is already defined over $L$, and $\sigma$ factors through the natural homomorphism  $\Gal(K^\sep/K) \to \Gal(L^\sep/L)$, where $L^\sep$ is the separable closure of $L$ in~$K^\sep$. Since $K$ is finitely generated over~$L$, the intersection $K\cap L^\sep$ is finite over~$L$; hence the image of this homomorphism is open. To prove the theorem we may thus replace $K$ by~$L$, after which $K$ is the residue field of~$\eta$.

\medskip
The assumption on $\psi$ implies that even after this reduction, $K$ is not a finite extension of~$F$. Therefore its transcendence degree over $F$ is $\ge1$. Let $k$ denote the algebraic closure of $F$ in~$K$. Then $\eta$ can be viewed as the generic point of a geometrically irreducible and reduced locally closed algebraic subvariety $X\subset M_k$ of dimension $\ge1$. After shrinking $X$ we may assume that $X$ is smooth. We are then precisely in the situation of the preceding section, with $\psi=\phi_\eta$. 
The homomorphism $\sigma$ above is then the composite
$$\Gal(\Ksep/K) \cong \pi_1(\eta,\bar\eta) 
\onto \pi_1(X,\bar\eta) \stackrel\rho\to \GL_r(\BA_F^f)$$
with $\rho$ as in Section~\ref{2}. It follows that the groups called $\Gamma$ in this section and the last coincide.
The desired openness is now equivalent to Theorem \ref{Main1}~(b).
\end{Proof}

{\bf Note:} The ad\`elic openness for a Drinfeld module defined over a {\em finite} 
extension of $F$ is still unproved.

\section{$\Fp$-Adic openness}
\label{p}

This section and the next three are devoted to proving Theorem~\ref{AnZD}.
Throughout we retain the notations from Sections \ref{1} and~\ref{2}
and the assumptions $\dim X\ge1$
and $\End_{\bar\eta_\BC}(\phi_{\bar\eta_\BC}) = A$.
In this section we recall a known result on $\Fp$-adic openness.
For any place $\Fp\not=\infty$ of $F$ let $\Gamma_\Fp$ denote the image of $\Gamma$ under the projection $\GL_r(\BA_F^f) \onto \GL_r(F_\Fp)$.

\begin{Thm}
\label{DMT1}
$\Gamma_\Fp$ is open in $\GL_r(F_\Fp)$.
\end{Thm}

\begin{Proof}
By construction $\Gamma_\Fp$ is the image of the monodromy representation
$$\rho_\Fp\!: \pi_1(X,\bar\eta) \longto \GL_r(F_\Fp)$$
on the rational $\Fp$-adic Tate module of~$\phi_{\eta}$. This is the same as the image of the composite homomorphism
$$\Gal(\Ksep/K) \cong \pi_1(\eta,\bar\eta) 
\onto \pi_1(X,\bar\eta) \stackrel{\rho_\Fp}\to \GL_r(F_\Fp).$$
Since $K$ is a finitely generated extension of~$F$,
and $\End_{K^\sep}(\phi_{\eta}) = A$ by the assumption and Lemma \ref{EE},
the desired openness is a special case of \cite[Thm.$\,$0.1]{PinkDMT}.
\end{Proof}


Next let $\Gamma^\geom_\Fp$ denote the image of $\Gamma^\geom$ under the projection $\GL_r(\BA_F^f) \onto \GL_r(F_\Fp)$. Note that this is a normal subgroup of~$\Gamma_\Fp$.
Lemma \ref{claim0a} immediately implies:

\begin{Lem}
\label{claim0b}
$\Gamma^\geom_\Fp$ is the closure of $g^{-1}\Delta_\Xi\, g$ in $\SL_r(F_\Fp)$
for some element $g \in \GL_r(F_\Fp)$.
\end{Lem}


\section{Zariski density}
\label{a}



\begin{Lem}
\label{claim1}
The Zariski closure $H$ of $\Delta_\Xi$ in $\GL_{r,F}$
is a normal subgroup of~$\GL_{r,F}$.
\end{Lem}

\begin{Proof}
Choose a place $\Fp\not=\infty$ of~$F$.
Then by base extension $H_{F_\Fp}$ is the Zariski closure of $\Delta_\Xi$ in $\GL_{r,F_\Fp}$.
Thus Lemma \ref{claim0b} implies that 
$g^{-1}H_{F_\Fp}g$ is the Zariski closure of $\Gamma^\geom_\Fp$ in $\GL_{r,F_\Fp}$.
Since $\Gamma_\Fp$ normalizes $\Gamma^\geom_\Fp$, it therefore normalizes $g^{-1}H_{F_\Fp}g$.
But $\Gamma_\Fp$ is open in $\GL_r(F_\Fp)$ by Theorem \ref{DMT1} and therefore Zariski dense in $\GL_{r,F_\Fp}$. Thus $\GL_{r,F_\Fp}$ normalizes $g^{-1}H_{F_\Fp}g$ and hence $H_{F_\Fp}$, and the result follows.
\end{Proof}


\begin{Lem}
\label{claim2}
$\Delta_\Xi$ is infinite.
\end{Lem}

\begin{Proof}
Let $X,K,k$ and $\phi_\eta$ be as in Section \ref{2}.
Then, as $M_k$ is affine and $\dim X\ge1$, 
there exists a valuation $v$ of $K$, corresponding to a point on the boundary of $X$
not on $M_k$, at which $\phi_\eta$ does not have potential good reduction. Denote by 
$I_v\subset\Gal(K^\sep/K k^\sep)$ the inertia group at~$v$. By the criterion of
N\'eron-Ogg-Shafarevich \cite[\S4.10]{GossBook}, the image of $I_v$ 
in $\Gamma^\geom_\Fp$ is infinite for any place $\Fp\not=\infty$ of~$F$.
In particular, $\Delta_\Xi$ is infinite by Lemma \ref{claim0b}, as desired.

\medskip
Alternatively, we may argue as follows. Suppose that $\Delta_\Xi$ is finite.
Then after increasing the
level structure we may assume that $\Delta_\Xi=1$.
Then $\Gamma^\geom_\Fp=1$ by Lemma \ref{claim0b}, which means that $\rho_\Fp$ factors as
$$\pi_1(X,\bar\eta) \ontoover{\ } \Gal(k^\sep/k) \longto \GL_r(F_\Fp).$$
After a suitable finite extension of the constant field $k$ we may assume that $X$ possesses a $k$-rational point~$x$.
Let $\phi_x$ denote the Drinfeld module over $k$ corresponding to~$x$.
Via the embedding $k\subset K$ we may consider it as a Drinfeld module over $K$ and compare it with~$\phi_\eta$.
The factorization above implies that the Galois representations on the $\Fp$-adic Tate modules of $\phi_x$ and $\phi_\eta$ are isomorphic.
By the Tate conjecture (see \cite{TaguchiTate} or \cite{Tamagawa2}) 
this implies that there exists an isogeny $\phi_x\to \phi_\eta$ over~$K$. Its kernel is finite and therefore defined over some finite extension $k'$ of~$k$. Thus $\phi_\eta$, as a quotient of $\phi_x$ by this kernel, is isomorphic to a Drinfeld module defined over~$k'$. But the assumption $\dim X\ge1$ implies that $\eta$ is not a closed point of~$M_k$; hence $\phi_\eta$ cannot be defined over a finite extension of~$k$. This is a contradiction.
\end{Proof}



\begin{Prop}
\label{AnZDa}
$\Delta_\Xi$ is Zariski dense in $\SL_{r,F}$.
\end{Prop}

\begin{Proof}
By construction we have $H\subset\SL_{r,F}$, and Lemma \ref{claim2} implies that $H$ is not contained in the center of $\SL_{r,F}$.
From Lemma \ref{claim1} it now follows that $H=\SL_{r,F}$, as desired.
\end{Proof}


The above results may be viewed as 
analogues of Andr\'e's results \cite[Thm.$\,$1, Prop.$\,$2]{Andre}, comparing the monodromy
group of a variation of Hodge structures with its generic Mumford-Tate group.
Our analogue of the former is $\Delta_\Xi$, and by \cite{PinkDMT} the latter corresponds to~$\GL_{r,F}$.
In our situation, however, we do not need the existence of a special point on $X$.



\section{Fields of coefficients}
\label{b}

Let $\bar\Delta_\Xi$ denote the image of $\Delta_\Xi$ in $\PGL_r(F)$. In this section we show that the field of coefficients of $\bar\Delta_\Xi$ cannot be reduced. 

\begin{Def}
\label{RC}
Let $L_1$ be a subfield of a field~$L$. We say that a subgroup $\bar\Delta \subset \PGL_r(L)$ lies in a model of $\PGL_{r,L}$ over $L_1$, if there exist a linear algebraic group $G_1$ over $L_1$ and an isomorphism  $\lambda_1\!: G_{1,L} \isoto \PGL_{r,L}$, such that $\bar\Delta \subset \lambda_1(G_1(L_1))$.
\end{Def}

\begin{Prop}
\label{AnZDb}
$\bar\Delta_\Xi$ does not lie in a model of $\PGL_{r,F}$ over a proper subfield of~$F$.
\end{Prop}

\begin{Proof}
As before we use an arbitrary auxiliary place $\Fp\not=\infty$ of~$F$. Let $\bar\Gamma^\geom_\Fp \triangleleft \bar\Gamma_\Fp$ denote the images of $\Gamma^\geom_\Fp \triangleleft \Gamma_\Fp$ in $\PGL_r(F_\Fp)$. Lemma \ref{claim0b} implies that $\bar\Gamma^\geom_\Fp$ is conjugate to the closure of $\bar\Delta_\Xi$ in $\PGL_r(F_\Fp)$. By Proposition \ref{AnZDa} it is therefore Zariski dense in $\PGL_{r,F_\Fp}$. On the other hand Theorem \ref{DMT1} implies that $\bar\Gamma_\Fp$ is an open subgroup of $\PGL_r(F_\Fp)$. It therefore does not lie in a model of $\PGL_{r,F_\Fp}$ over a proper subfield of~$F_\Fp$. Thus $\bar\Gamma^\geom_\Fp$ is Zariski dense and normal in a subgroup that does not lie in a model over a proper subfield of~$F_\Fp$, which by \cite[Cor.$\,$3.8]{PinkCS} implies that $\bar\Gamma^\geom_\Fp$, too, does not lie in a model over a proper subfield of~$F_\Fp$.

\medskip
Suppose now that $\bar\Delta_\Xi \subset \lambda_1(G_1(F_1))$ for a subfield $F_1\subset F$, a linear algebraic group $G_1$ over $F_1$, and an isomorphism  $\lambda_1\!: G_{1,F} \isoto \PGL_{r,F}$. Since $\bar\Delta_\Xi$ is Zariski dense in $\PGL_{r,F}$, it is in particular infinite. Therefore $F_1$ must be infinite. As $F$ is finitely generated of trans\-cen\-dence degree $1$ over~$\BF_p$, it follows that $F_1$ contains a transcendental element, and so $F$ is a finite extension of~$F_1$. Let $\Fp_1$ denote the place of $F_1$ below~$\Fp$. Since $\bar\Gamma^\geom_\Fp$ is the closure of $\bar\Delta_\Xi$ in $\PGL_r(F_\Fp)$, it is contained in $\lambda_1(G_1(F_{1,\Fp_1}))$. The fact that $\bar\Gamma^\geom_\Fp$ does not lie in a model over a proper subfield of~$F_\Fp$ thus implies that $F_{1,\Fp_1} = F_\Fp$.

\medskip
But for any proper subfield $F_1\subsetneqq F$, we can choose a place $\Fp\not=\infty$ of $F$ above a place $\Fp_1$ of~$F_1$, such that the local field extension $F_{1,\Fp_1} \subset F_\Fp$ is non-trivial. Thus we must have $F_1=F$, as desired.
\end{Proof}


\section{Strong approximation}
\label{c}


The remaining ingredient is the following general theorem.

\begin{Thm}
\label{SA}
For $r\ge2$ let $\Delta\subset \SL_r(F)$ be a subgroup that is contained in a congruence subgroup commensurable with $\SL_r(A)$. Assume that $\Delta$ is Zariski dense in $\SL_{r,F}$ and that its image $\bar\Delta$ in $\PGL_r(F)$ does not lie in a model of $\PGL_{r,F}$ over a proper subfield of~$F$. Then the closure of  $\Delta$  in  $\SL_r(\BA_F^f)$  is open.
\end{Thm}

\begin{Proof}
For finitely generated subgroups this is a special case of \cite[Thm.$\,$0.2]{PinkSA}. 
That result concerns arbitrary finitely generated Zariski dense subgroups of $G(F)$ for arbitrary semisimple algebraic groups~$G$, but it uses the finite generation only to guarantee that the subgroup is integral at almost all places of~$F$. For $\Delta$ as above the integrality at all places $\not=\infty$ is already known in advance, so the proof in \cite{PinkSA} covers this case as well.

\medskip
As an alternative, we will deduce the general case by showing that every sufficiently large finitely generated subgroup $\Delta_1 \subset \Delta$ satisfies the same assumptions. Then the closure of $\Delta_1$ in $\SL_r(\BA_F^f)$ is open by \cite{PinkSA}, and so the same follows for~$\Delta$, as desired.

\medskip
For the Zariski density of $\Delta_1$ note first that the trace of the adjoint representation defines a dominant morphism to the affine line $\SL_{r,F} \to \BA^1_F$,  $g\mapsto \trace(\Ad(g))$.
Since $\Delta$ is Zariski dense, this function takes infinitely many values on~$\Delta$. As the field of constants in $F$ is finite, we may therefore choose an element $\gamma \in \Delta$ with $\trace(\Ad(\gamma))$ transcendental. Then $\gamma$ has infinite order; hence the Zariski closure $H\subset\SL_{r,F}$ of the abstract subgroup generated by $\gamma$ has positive dimension. Let $H^\circ$ denote its identity component. Since $\Delta$ is Zariski dense and $\SL_{r,F}$ is almost simple, the $\Delta$-conjugates of $H^\circ$ generate $\SL_{r,F}$ as an algebraic group. By noetherian induction finitely many conjugates suffice. It follows that finitely many conjugates of $\gamma$ generate a Zariski dense subgroup of $\SL_{r,F}$. Thus every sufficiently large finitely generated subgroup $\Delta_1 \subset \Delta$ is Zariski dense.

\medskip
Consider such $\Delta_1$ and let $\bar\Delta_1$ denote its image in $\PGL_r(F)$. Consider all triples $(F_1,G_1,\lambda_1)$ consisting of a subfield $F_1\subset F$, a linear algebraic group $G_1$ over $F_1$, and an isomorphism  $\lambda_1\!: G_{1,F}\isoto \PGL_{r,F}$, such that $\bar\Delta_1 \subset \lambda_1(G_1(F_1))$. By \cite[Thm.$\,$3.6]{PinkCS} there exists such a triple with $F_1$ minimal, and this $F_1$ is unique, and $G_1$ and $\lambda_1$ are determined up to unique isomorphism. Consider another finitely generated subgroup  $\Delta_1 \subset \Delta_2 \subset \Delta$ and let $(F_2,H_2,\lambda_2)$ be the minimal triple associated to it. Then the uniqueness of $(F_1,G_1,\lambda_1)$ implies that $F_1\subset F_2$, that $G_2 \cong G_{1,F_2}$, and that $\lambda_2$ coincides with the isomorphism $G_{2,F} \cong G_{1,F} \to \PGL_{r,F}$ obtained from~$\lambda_1$. In other words, the minimal model $(F_1,G_1,\lambda_1)$ is monotone in~$\Delta_1$.

\medskip
For any increasing sequence of Zariski dense finitely generated subgroups of~$\Delta$ we thus obtain an increasing sequence of subfields of~$F$. This sequence must become constant, say equal to $F_1\subset F$, and the associated model of $\PGL_{r,F}$ over $F_1$ is the same up to isomorphism from that point onwards. Thus we have a triple $(F_1,G_1,\lambda_1)$ with $\bar\Delta_1 \subset \lambda_1(G_1(F_1))$ for every sufficiently large finitely generated subgroup $\bar\Delta_1 \subset \bar\Delta$. But then we also have $\bar\Delta \subset \lambda_1(G_1(F_1))$, which by assumption implies that $F_1=F$. Thus every sufficiently large finitely generated subgroup of $\Delta$ satisfies the same assumptions as $\Delta$, as desired.
\end{Proof}

\begin{Proofof}{Theorem \ref{AnZD}}
In the situation of Theorem \ref{AnZD} we automatically have $r\ge2$,
so the assertion follows by combining Propositions \ref{AnZDa} and \ref{AnZDb} with Theorem~\ref{SA} for $\Delta_\Xi$.
\end{Proofof}


\section{Arbitrary endomorphism rings}
\label{End}

Set $E := \End_{\bar\eta_\BC}(\phi_{\bar\eta_\BC})$, which is a finite integral ring extension of~$A$. 
Write $r=r'\cdot[E/A]$; then the centralizer of $E$ 
in $\GL_r(\BA_F^f)$ is isomorphic to $\GL_{r'}(E\otimes_A\BA_F^f)$.
Lemma \ref{EE} implies that all elements of $E$ are defined over some fixed finite extension of~$K$. This means that an open subgroup of $\rho\bigl( \pi_1(X,\bar\eta) \bigr)$
is contained in $\GL_{r'}(E\otimes_A\BA_F^f)$. 
Thus by Lemma \ref{claim0a} the same holds for a subgroup of finite index of $\Delta_\Xi$.
The following results can be deduced easily from Theorems \ref{AnZD}, \ref{Main1}, and~\ref{DMT2}, using the same arguments as in \cite[end of \S2]{PinkDMT}.

\begin{Thm}
\label{AnZDE}
In the situation of before Theorem \ref{AnZD}, for $E := \End_{\bar\eta_\BC}(\phi_{\bar\eta_\BC})$ arbitrary,
the closure in $\GL_r(\BA_F^f)$ of some subgroup of finite index of $\Delta_\Xi$
is an open subgroup of $\SL_{r'}(E\otimes_A\BA_F^f)$.
\end{Thm}

\begin{Thm}
\label{Main2}
In the situation of before Theorem \ref{Main1}, for $E := \End_{K^\sep}(\phi_\eta)$ arbitrary,
\begin{itemize}
\item[(a)] some open subgroup of $\Gamma^\geom := \rho\bigl( \pi_1(X_{k^\sep},\bar\eta) \bigr)$ is an open subgroup of $\SL_{r'}(E\otimes_A\BA_F^f)$, and
\item[(b)] some open subgroup of $\Gamma := \rho\bigl( \pi_1(X,\bar\eta) \bigr)$ is an open subgroup of $\GL_{r'}(E\otimes_A\BA_F^f)$.
\end{itemize}
\end{Thm}

\begin{Thm}
\label{DMT2E}
In the situation of before Theorem \ref{DMT2}, for $E := \End_{K^\sep}(\psi)$ arbitrary,
suppose that $\psi$ cannot be defined over a finite extension of $F$ inside $K^\sep$.
Then some open subgroup of $\Gamma := \sigma\bigl( \Gal(K^\sep/K) \bigr)$ is an open subgroup of $\GL_{r'}(E\otimes_A\BA_F^f)$.
\end{Thm}



\end{document}